\documentclass[a4paper, 10pt, conference]{ieeeconf}    

\IEEEoverridecommandlockouts  
\usepackage[a4paper, left=1.91cm, right=1.31cm, top=1.91cm, bottom=3.67cm]{geometry}
\usepackage[utf8]{inputenc}
\usepackage{amsmath,amsfonts,amssymb}
\usepackage{mathtools}
\usepackage{dsfont}
\usepackage[svgnames]{xcolor}
\usepackage{graphicx}
\usepackage{float}
\usepackage{caption}
\usepackage{subcaption}
\usepackage{tikz}
\usepackage[expansion=false]{microtype}
\usepackage{url}
\usepackage{comment}

\newcommand{\nl}{\left\|}
\newcommand{\nr}{\right\|}

\newcommand{\Norm}[2][ ]{\nl #2 \nr_{#1}}

\newcommand{\N}{\mathds{N}}

\newcommand{\R}{\mathds{R}}

\renewcommand{\phi}{\varphi}

\makeatletter
\newcommand{\Itemlabel}[2]{#2\def\@currentlabel{#2}\label{#1}}
\makeatother

\makeatletter
\newcommand\setcurrentname[1]{\def\@currentlabelname{#1}}
\makeatother

\newcommand{\nocontentsline}[3]{}
\newcommand{\tocless}[2]{\bgroup\let\addcontentsline=\nocontentsline#1{#2}\egroup}
\newcommand{\toclesslab}[3]{\bgroup\let\addcontentsline=\nocontentsline#1{#2\label{#3}}\egroup} 

\newtheorem{thm}{Theorem}

\newtheorem{definition}[thm]{Definition}
\newtheorem{remark}[thm]{Remark}

\title{Event-Triggered Funnel Control for Output Reference Tracking of Unknown Nonlinear Systems }

\author{Janina Schaa$^{1}$ and Thomas Berger$^{1}$%
\thanks{*This work was funded by the Deutsche Forschungsgemeinschaft (DFG, German Research Foundation)~-- Project-ID 524064985.}
\thanks{$^{1}$Institut f\"ur Mathematik, Martin-Luther-Universit\"at Halle-Wittenberg, Halle (Saale), Germany
        {\tt\small \{janina.schaa, thomas.berger\}@mathematik.uni-halle.de}}%
}

\begin{document}
\bstctlcite{IEEEexample:BSTcontrol}
\maketitle
\thispagestyle{empty}
\pagestyle{empty}

\begin{abstract}
We propose an event-triggered zero-order hold implementation of the funnel controller for output reference tracking of unknown nonlinear systems. The control input is updated only when the deviation between the current and the last transmitted input exceeds a prescribed threshold, reducing the number of input updates while maintaining closed-loop performance guarantees by incorporating an auxiliary integral state to compensate for the error introduced by the zero-order-hold input. We establish a uniform lower bound on the inter-event intervals excluding Zeno behavior, and derive a growth condition on the nonlinear dynamics that ensures global existence of the solution for arbitrary, predefined sampling schedules. The theoretical results are illustrated by numerical simulations.
\end{abstract}

\section{Introduction}

Reference tracking is a fundamental control objective arising throughout engineering applications, and in safety-critical situations it is not sufficient to guarantee asymptotic convergence, but the tracking error must remain within prescribed bounds for all times. Funnel control, introduced in~\cite{Ilchmann2002}, is an adaptive high-gain control framework that achieves this objective under mild structural assumptions, see~\cite{SurveyBerger}. The practical applicability of funnel controllers to real-world problems has been demonstrated in several engineering domains, including the control of industrial servo systems~\cite{Hackl2017,Ilchmann2009}, DC-link power flow control~\cite{Senfelds2014} and voltage and current control of electrical circuits~\cite{Berger2014}.

While related approaches such as prescribed performance control~\cite{Bechlioulis2008,bechlioulis_2014} utilize a nonlinear transformation of the tracking error and design a controller based on the transformed error, funnel controllers incorporate the performance constraint directly into the design of the feedback law.

Several challenges in the application of funnel controllers have been addressed in recent research, including output measurement losses~\cite{Berger2023}, input constraints~\cite{Hu2022,Berger2024}, discrete-time systems~\cite{Cheng2023}, and the application to systems of higher relative degree~\cite{BergerReis18}. A practically relevant limitation that has received less attention is the requirement of a continuously varying control input. Classical funnel controllers rely on continuous-time actuation in order to enforce the funnel constraint, which is difficult to implement in practice. Maintaining the performance constraint despite infrequent input updates and without knowledge of the system dynamics is therefore a nontrivial challenge.

In the funnel control literature, a controller with zero-order hold structure has been proposed in~\cite{lanza2024b,lanza2024a}, utilizing knowledge of upper bounds on the disturbance and the dynamics of the system in the controller design. In the prescribed performance control framework, sampled-data implementations have been studied via data-driven approaches approximating the unknown nonlinear dynamics~\cite{liu2024,li2018}. Moreover, event-triggered prescribed performance control has been applied to consensus problems in multi-agent systems~\cite{yao2025}. 

Event-triggered control has emerged as a practically motivated paradigm for reducing the communication and computational burden associated with continuous-time systems~\cite{ETCSurvey}. Rather than updating the control input at fixed time instants as in periodic sampling, event-triggered schemes determine the next update time online based on a state- or output-dependent trigger condition, generally allowing for sparser sampling, see~\cite{Peng2018} and the references therein. A central requirement in any event-triggered scheme is the exclusion of Zeno behavior, the accumulation of infinitely many events in finite time~\cite{GoebelHybridSystems}.

In this paper, we address the limitation of infrequent controller updates by proposing an event-triggered implementation of the funnel controller, in which the control input is updated whenever the deviation between the current and the last transmitted input exceeds a threshold $M>0$. This event-triggering mechanism is particularly well-suited to the model-free nature of funnel control. Since no system knowledge is available, the controller cannot anticipate a proper sampling schedule. The event-triggered scheme automatically adapts the sampling rate to the current system behavior, triggering updates whenever needed to maintain the performance constraint at all times, and refraining from unnecessary updates otherwise. 

The proposed controller incorporates an auxiliary state variable $I$ that dynamically compensates for the error introduced by the zero-order hold input applied to the system, allowing the funnel constraint to be maintained despite infrequent updates. We establish a lower bound on the length of the inter-event intervals to exclude Zeno behavior, and further show that predefined sampling schedules yield global solutions, if the system dynamics satisfy a growth condition, connecting this approach to time-triggered implementations.

We use the following notation. $\N$ and $\R$ denote the natural and real numbers, respectively, $\N_0 = \N \cup \{ 0 \}$, and $\R_{\geq 0} = [0,\infty)$. $\Norm{x}$ denotes the Euclidean norm of $x\in\R^n$. $C(V,\R^n)$ is the linear space of continuous functions $f:V\rightarrow \R^n$, where $V \subseteq \R^m$. On an interval $I \subseteq \R$, $L^{\infty}(I,\R^n)$ denotes the space of measurable and essentially bounded functions $f:I\rightarrow \R^n$ with norm $\Vert f\Vert_{\infty} = \operatorname{ess} \sup_{t \in I} \Vert f(t) \Vert$, while $L_{\text{loc}}^{\infty}(I,\R^n)$ is the set of measurable and locally essentially bounded functions. Moreover, $W^{k,\infty}(I,\R^n)$ is the Sobolev space of all $k$-times weakly differentiable functions $f:I\rightarrow \R^n$ such that $f,\ldots,f^{(k)} \in L^{\infty}(I,\R^n)$. 

\section{Problem Formulation}

We consider nonlinear first order systems
\begin{equation} \label{eq:System}
\dot{y}(t) = f(y(t)) + S_Mu(t), \quad y(0) = y^0 \in \R^n,
\end{equation}
given by an unknown continuous function $f \in C(\R^n,\R^n)$ and an initial condition $y^0 \in \R^n$.  The objective is to ensure that the tracking error $y(t) - y_{\text{ref}}(t)$ for a given trajectory $y_{\text{ref}} \in W^{1,\infty}(\R_{\geq 0}, \R^n)$ remains within prescribed performance bounds for all times, where the control input can only be updated at discrete sampling instances rather than continuously. The continuous control input $u \in C(\R_{\geq 0} , \R^n)$ is applied via the event-triggered sampling mechanism described by the operator $S_M$ in~\eqref{eq:System}, defined in the following, for some design parameter $M > 0$. The sampling mechanism should further include a lower bound on the time between consecutive updates in order to exclude Zeno behavior.

In order to define $S_M$, let $u: [0,\tau) \rightarrow \R^n$ be a continuous signal for some $\tau \in (0,\infty]$. 
\begin{definition}
For a parameter $M > 0$, let the update map $\rho_M : \R^n \times \R^n \rightarrow \R^n$ be given by
\begin{equation*}
\rho_M(y,m) = \begin{cases}  m&, \quad \text{ if } \Norm{y-m} < M \\ y&, \quad \text{ if } \Norm{y-m} \geq M.  \end{cases}
\end{equation*}
Then, the event-triggered sampling operator $S_M: C([0,\tau), \R^n) \rightarrow L_{\text{loc}}^{\infty}([0,\tau), \R^n)$ is defined via $(S_Mu)(0) = u(0),$ and
\begin{equation}
(S_Mu)(t) = \rho_M(u(t), (S_Mu)(t_-)), \quad \forall t \in (0,\tau),
\end{equation}
where $(S_Mu)(t_-) = \lim_{s \uparrow t} (S_Mu)(s)$ denotes the left-hand limit.
\end{definition}

The sampling operator satisfies $\Vert y - \rho_M(y,m) \Vert < M$, and $\rho_M(y,y) = y$, for all $y,m \in \R^n$.

Let $(t_k)_{k \in T}$ denote the increasing sequence of points of discontinuity of $S_Mu$, where $T \in \left\{ \{1,\hdots,N \} \mid N \in \N \right\}  \cup \{ \N \} $, which are the update times generated by the map $\rho_M$, and set $t_0 = 0$. At every such time $t_k$, the update condition implies 
$$ \Vert u(t_k) - u(t_{k-1}) \Vert = \Vert u(t_k) - (S_Mu)(t_k^-) \Vert \geq M.$$
We claim that the sequence of update times is locally finite, that is, either $T$ is finite, or $\lim_{k \rightarrow \infty} t_k = \tau$. Fix an arbitrary compact subinterval $[0,t] \subseteq [0,\tau)$. Then $u$ is uniformly continuous on the compact interval $[0,t]$, hence $\Norm{u(t_k) - u(t_{k-1})} \geq M$ implies that $t_k -  t_{k-1} \geq \delta$ for some $\delta > 0$ independent of $k$. Therefore, only finitely many points of discontinuity of $S_Mu$ can exist in the interval $[0,t]$. In particular, there is no accumulation point of update times in $[0,\tau)$. Hence, $S_Mu$ is well-defined on $[0,\tau)$, as the left-hand limits $(S_Mu)(t_-)$ exist for all $t \in (0,\tau)$. We refer to the intervals $[t_{k-1},t_k)$ as sampling intervals.

\section{Controller Design}

Choose a dynamic control law
\begin{equation} \label{eq:Controller}
\boxed{
\begin{aligned}
u(t) &= -k(t,e(t))e(t),\\
k(t,e) &= \frac{1}{1 - \varphi(t)^2\Norm{e}^2}, \\
e(t) &= y(t) -y_{\text{ref}}(t) + I(t), \\
\dot{I}(t) &= -\alpha I(t) +u(t) -S_Mu(t) , &\ I(0) = I^0 \in \R^n,
\end{aligned}
}
\end{equation}
where $I$ can be interpreted as a dynamical adjustment of the reference trajectory $y_{\text{ref}} \in W^{1,\infty}(\R_{\geq 0}, \R^n)$, and $\varphi \in W^{1,\infty}(\R_{\geq 0}, \R)$ with $\inf_{t \geq 0} \varphi(t) > 0$. The constants $\alpha, M>0$ are design parameters of the controller. The control law can be interpreted as an event-triggered sample-and-hold implementation of a funnel controller.

A function $(y,I):[0,\omega) \rightarrow \R^n \times \R^n$, $\omega \in (0,\infty]$, is called solution of the closed-loop system $\eqref{eq:System},~\eqref{eq:Controller}$ in the sense of Carathéodory, if it satisfies the initial conditions $(y(0),I(0)) = (y^0,I^0)$, and $(y,I)$ is locally absolutely continuous on $[0,\omega)$ and satisfies
\begin{align*}
\dot{y}(t) &= f(y(t)) + S_Mu(t),\\
\dot{I}(t) &= -\alpha I(t) +u(t) -S_Mu(t), 
\end{align*}
with $u$ as in~\eqref{eq:Controller} for almost all $t \in [0,\omega)$. We call a solution maximal, if it does not have a right extension, which is also a solution. A global solution is one with $\omega=\infty$.

\begin{thm} \label{Thm:event}
Consider the closed-loop system~\eqref{eq:System},~\eqref{eq:Controller}. Assume that the initial values satisfy $\varphi(0)\Norm{y^0-y_{\text{ref}}(0) + I^0} < 1$. Then the closed-loop system has a global solution satisfying
\begin{enumerate}
\item $\varphi(t)\Norm{e(t)} \leq \varepsilon$ for all times $t \geq 0$ for some $\varepsilon \in (0,1)$ depending on $M$ and $\alpha$,
\item there exists a uniform lower bound $t_{\text{min}}>0$ depending on $M$ and $\alpha$ for the length of the intervals where $S_M u$ is constant, 
\item $u$ remains bounded, and
\item $\Vert I \Vert_{\infty} \leq \max \{ \frac{M}{\alpha}, \Norm{I^0} \}$. 
\end{enumerate}
\end{thm} 
\begin{proof}
\emph{Step 1:} We show the existence of a maximal solution $(y,I):[0,\omega) \rightarrow \R^n \times \R^n$.\\
We will inductively show the existence of a solution on the intervals on which $S_Mu$ is constant. Define the non-empty open set
$$\mathcal{D} \!\coloneq \!\left\{  (t,y,I)\! \in \!\R_{\geq 0}\! \times\! \R^n \!\times \!\R^n \mid \varphi(t)\Norm{y-y_{\text{ref}}(t) + I}\! <\! 1 \right\}\!. $$
Define the function 
\begin{align*}
G_v:\ &\mathcal{D} \rightarrow \R^n \times \R^n, \\
 (t,y,I) &\mapsto \begin{pmatrix} f(y)+v\\  -\alpha I+\frac{-(y-y_{\text{ref}}(t) +I)}{1-\varphi(t)^2\Norm{y-y_{\text{ref}}(t) +I}^2}-v\end{pmatrix},
\end{align*}
for some fixed $v \in \R^n$. 
Then, the initial value problem~\eqref{eq:System},~\eqref{eq:Controller} takes the form
\begin{equation} \label{eq:ClosedLoop}
\begin{pmatrix} \dot{y}(t) \\ \dot{I}(t) \end{pmatrix} = G_v(t,y(t),I(t)), \quad \begin{pmatrix} y(s) \\ I(s) \end{pmatrix} = \begin{pmatrix} y^s \\ I^s \end{pmatrix}
\end{equation}
on an interval $[s,\tau)$ such that there are no points of discontinuity of $S_Mu$ in this interval, with $v = S_Mu(s)$. We observe that $G_v$ is continuous in $(y,I)$, and measurable in $t$, and locally bounded on compact subsets of $\mathcal{D}$. Whenever $(s, y^{s}, I^{s})\in \mathcal{D}$, application of~\cite[\S\,10, Thm.~XX]{Walt98} yields the existence of a maximal solution $(y,I):[s,\tilde\omega) \rightarrow \R^n \times \R^n$ for some $\tilde\omega \in (s,\infty]$. Moreover, the closure of the graph of the maximal solution is not contained in a compact subset of $\mathcal{D}$.

If $\Norm{u(t) - v} < M$ for all $t \in [s,\tilde\omega)$ where $u$ is defined as in~\eqref{eq:Controller}, then $\frac{\Norm{e(t)}}{1-\varphi(t)^2\Norm{e(t)}^2} \leq \Norm{v} + M$ for all $t \in [s,\tilde\omega)$. This implies that $u(\cdot)$ is uniformly bounded on this interval. Hence, either $\tilde{\omega} = \infty$, or the solution is contained in a compact subset of $\mathcal{D}$, which is a contradiction. This shows that either $\tilde\omega = \infty$, or there exists $t \in [s,\tilde\omega)$ such that $\Norm{v - u(t)} \geq M$.

We inductively define the sequence $(t_k)_{k \in T}$ for some index set $T \subseteq \N$ as follows. Set $t_0 = 0$, and observe that $(t_0,y^0,I^0) \in \mathcal{D}$. For each $k \in \N_0$, apply the above argument for the existence of a maximal solution on some interval $[t_k,\tilde\omega)$. If $\Norm{u(t) - u(t_k)} < M$ for all $t \in [t_k,\tilde\omega) $, then $\tilde\omega = \infty$ by the previous argument, and the induction terminates with $T = \{ 1, \ldots, k\}$. Otherwise, choose $t_{k+1} \coloneq \inf \{ t\in [t_k,\tilde\omega) \mid \Norm{u(t) - u(t_k)} = M \}$, which is well-defined by the previous argument and the continuity of $u$. Moreover, we observe that $(t_{k+1},y(t_{k+1}),I(t_{k+1})) \in \mathcal{D}$ since $t_{k+1} < \tilde\omega$, closing the induction step. Set $T = \N$ if the induction did not terminate at any step $k \in\N$.

By construction, the sequence $(t_k)_{k \in T}$ is the sequence of discontinuities of $S_Mu$, and the closed-loop system~\eqref{eq:System},~\eqref{eq:Controller} coincides with~\eqref{eq:ClosedLoop} on the intervals $[t_k,t_{k+1})$ with $v=u(t_k)$. 

Hence, we can inductively conclude that there exists a maximal solution of the closed-loop system~\eqref{eq:System},~\eqref{eq:Controller} on $[0,\omega)$ for $\omega = \lim_{k \rightarrow \infty} t_k$, if $T = \N$, and $\omega = \infty$ if $\vert T \vert < \infty$. Note that, at this point, $\omega< \infty$ is not excluded in the case $T=\N$.

\emph{Step 2:} We show boundedness of the displacement $I$.\\  
Aiming for a contradiction, assume that there is $\tau_0 \in (0,\omega)$ such that $\Norm{I(\tau_0)} > \max \{ \frac{M}{\alpha}, \Norm{I^0} \} =: \tilde{M}$. Choose $\tau^{\star} \coloneq \sup \{ t \leq \tau_0 \mid \Norm{I(t)} = \tilde{M} \}, $ which is well-defined by continuity of $I$ and the choice of $\tilde{M}$. 
Since by construction of the operator $S_M$, the estimate
\begin{equation} \label{eq:SamplingCondition}
\Norm{u(t) - S_Mu(t)} \leq M
\end{equation}
holds, we can conclude
\begin{align*}
\frac{1}{2} \frac{d}{dt} \Norm{I(t)}^2 &= \langle I(t), -\alpha I(t)+u(t)-S_Mu(t) \rangle\\
&\leq \Norm{I(t)} (- \alpha \Norm{I(t)} + \Norm{u(t) - S_Mu(t)} )\\
&\leq  \Norm{I(t)} (-M+M) = 0 
\end{align*}
 for all $t \in [\tau^{\star},\tau_0]$. Integration yields
$$ \tilde{M} = \Norm{I(\tau^{\star})} \geq \Norm{I(\tau_0)} > \tilde{M}, $$
a contradiction. Hence, we can conclude that $\Norm{I(t)} \leq \tilde{M}$ holds for all $t \in [0,\omega)$. 

\emph{Step 3:} We show $\varphi(t) \Norm{e(t)} \leq \varepsilon$ and boundedness of $u$.
Due to its continuity, $f$ is uniformly bounded on each compact subset of $\R^n$. We define $F:\R_{\geq 0} \rightarrow \R_{\geq 0}, r \mapsto \left \Vert f_{\vert B_{r+\left\Vert \frac{1}{\varphi} \right\Vert_{\infty} + \Vert y_{\text{ref}} \Vert_{\infty}}(0)} \right\Vert_{\infty}$, where $B_r(0) \subseteq \R^n$ denotes the ball of radius $r$ centered at the origin. Then by construction, the estimate
\begin{equation} \label{eq:fEstimate}
\Vert f(y) \Vert \leq F(\tilde{M}), \quad \text{ for all } \Norm{y} \leq \Vert y_{\text{ref}} \Vert_{\infty} + \tilde{M} + \left\Vert \frac{1}{\varphi} \right\Vert_{\infty},
\end{equation}
holds, where $\tilde{M}$ is defined in the previous step. Choose the constant $\varepsilon \in (0,1) $ depending on the system parameters such that $\varphi(0)\Norm{y^0-y_{\text{ref}}(0) + I^0} < \varepsilon$ and
\begin{equation} \label{eq:eps}
\frac{\varepsilon^2}{1-\varepsilon^2} \geq  \left\Vert \frac{\dot{\varphi}}{\varphi} \right\Vert_{\infty} + \Vert \varphi \Vert_{\infty} \left( \Vert \dot{y}_{\text{ref}} \Vert_{\infty} + F(\tilde{M}) + \alpha\tilde{M} \right) .
\end{equation}

Assume that there exists $s_1 \in [0,\omega)$ such that $\varphi(s_1) \Norm{e(s_1)} > \varepsilon$. Then $s_0 \coloneq \sup \{ s \in [0,s_1] \mid \varphi(s)\Norm{e(s)} = \varepsilon \}$ is well-defined by the continuity of $e$ and $\varphi$, as well as the choice of $\varepsilon$. Then on $[s_0,s_1]$, the estimate
\begin{equation}\label{eq:relError}
\begin{aligned}  
\frac{1}{2} \frac{d}{dt} \varphi(t)^2 \Norm{e(t)}^2 &= \dot{\varphi}(t) \varphi(t) \Norm{e(t)}^2 + \varphi(t)^2 \langle e(t), \dot{e}(t)\rangle \\
&\overset{\eqref{eq:fEstimate}}{\leq} \Vert \varphi \Vert_{\infty} \left( \Vert \dot{y}_{\text{ref}} \Vert_{\infty} + F(\tilde{M}) + \alpha\tilde{M} \right) \\
&\quad \quad + \left\Vert \frac{\dot{\varphi}}{\varphi} \right\Vert_{\infty}  - \frac{\varepsilon^2}{1- \varepsilon^2} \leq 0
\end{aligned}
\end{equation}
holds. Integration yields
$$ \varepsilon = \varphi(s_0)\Norm{e(s_0)} \geq \varphi(s_1)\Norm{e(s_1)} >  \varepsilon,$$
a contradiction. We can conclude that $\varphi(s)\Norm{e(s)} \leq \varepsilon$ for all $s \in [0,\omega)$. Moreover, this yields $\Norm{u(t)} \leq \frac{1}{1-\varepsilon^2} \left\Vert \frac{1}{\varphi} \right\Vert_{\infty}$, hence $u$ remains bounded.

\emph{Step 4:} We derive a uniform lower bound on the length of the sampling intervals.\\
We can observe that, using $\Norm{I(t)} \leq \tilde{M}$ on $[0,\omega)$ from Step 2 and $\Norm{u(t)} \leq \frac{1}{1-\varepsilon^2} \left\Vert \frac{1}{\varphi} \right\Vert_{\infty}$ from Step 3, we can obtain the estimate
\begin{equation} \label{eq:Derivative}
\begin{aligned} 
\Norm{\frac{d}{dt} u(t)}\! &\leq \!\Norm{ \frac{\dot{e}(t)}{1\!-\!\varphi(t)^2 \Norm{e(t)}^2} } \!+ \!\Norm{\frac{\frac{d}{dt} (\varphi(t)^2 \Norm{e(t)}^2)}{(1\!-\!\varphi(t)^2\Norm{e(t)}^2)^2} e} \\
&\overset{\eqref{eq:relError}}{\leq} \frac{1}{1-\varepsilon^2} \left( \Vert \dot{y}_{\text{ref}} \Vert_{\infty} + F(\tilde{M}) + \alpha \tilde{M} + \frac{\Vert \frac{1}{\varphi} \Vert_{\infty}}{1-\varepsilon^2} \right) \\
&\quad +\! \frac{2 \Vert \frac{1}{\varphi} \Vert_{\infty}}{(1-\varepsilon^2)^2} \Bigg(\!\Vert \varphi \Vert_{\infty} \!\left( \Vert \dot{y}_{\text{ref}} \Vert_{\infty} \!+\! F(\tilde{M}) + \alpha \tilde{M} \right) \\
& \quad \quad \quad + \left\Vert \frac{\dot{\varphi}}{\varphi} \right\Vert_{\infty} + \frac{\varepsilon^2}{1- \varepsilon^2}\Bigg) =: c(\alpha, \tilde{M}). 
\end{aligned}
\end{equation}

This yields
\begin{equation*}
\begin{aligned}
\Norm{u(t_{k+1}) - u(t_k)} &= \Norm{\int_{t_k}^{t_{k+1}} \frac{d}{ds} u(s) ds} \\
&\leq (t_{k+1}-t_k) \Norm{\frac{d}{ds} u(s)} \\
&\leq (t_{k+1}-t_k) c(\alpha,\tilde{M}), 
\end{aligned}
\end{equation*}
for each $k \in \N_0$, hence 
\begin{equation*}
t_{k+1}-t_k \geq \frac{\Norm{u(t_{k+1}) - u(t_k)}}{c(\alpha, \tilde{M})} = \frac{M}{c(\alpha,\tilde{M})}.
\end{equation*}
In particular, this estimate ensures $t_k \rightarrow \infty$ as $k \rightarrow \infty$.

\emph{Step 5:} We show the existence of a global solution.\\
Assume that $\omega < \infty$. Then $T = \N$, and by Step~1, $\omega = \lim_{k \rightarrow \infty} t_k = \infty$, due to the uniform lower bound on the interval lengths $t_{k+1} - t_k$ established in Step~4. We can conclude $\omega = \infty$, finishing the proof.
\end{proof}

\begin{remark}
The presence of the correction term $I$ allows the controller to accommodate arbitrarily large initial tracking errors $y^0 - y_{\text{ref}}(0)$, provided that the initial tracking error is known with sufficient accuracy, in order to allow for a suitable $I^0$ that ensures $\Norm{e(0)} \varphi(0) < 1$.
\end{remark}

\begin{remark}
While the bound prescribed by $\varphi$ holds for $e(t) = y(t) - y_{\text{ref}}(t) + I(t)$, we can use the boundedness of $I$ to conclude that the tracking error satisfies $\Vert y(t) - y_{\text{ref}}(t) \Vert < \varphi(t)^{-1} + \tilde M$. Since $\tilde M>0$ is determined through design parameters $M,\alpha, \Norm{I^0}$, the right hand side can resemble an arbitrary bound $\psi \in W^{1,\infty}(\R_{\geq 0}, \R_{\geq 0})$ with $\inf_{t \geq 0} \psi(t) > 0$.
\end{remark}

\begin{remark}
The choice of $\alpha$ and $M$ entails a trade-off between the admissible sampling interval lengths and the resulting transient performance bounds. Choosing the parameters such that $\frac{M}{c(\alpha,\tilde{M})}$ is maximized yields the least restrictive a priori bound on the sampling interval length and hence allows for the coarsest sampling scheme. On the other hand, choosing $\tilde{M}$ small reduces the offset from the prescribed performance bound, since $\Vert I \Vert_{\infty} \leq \tilde{M}$. Thus, when tighter safety margins are desired, $\tilde{M}$ can be decreased, at the expense of a more restrictive bound on the lengths of the sampling intervals.
\end{remark}

\section{Predefined Sampling Schedules}

Instead of monitoring the control signal continuously online to determine the next sampling time, one may alternatively fix the sampling sequence $(t_k)_{k\in\mathbb N_0}$ prior to implementation. In contrast to the event-triggered controller, where the sampling times are determined online by the input trajectory $u(t)$, the time-triggered sampling sequence is chosen a priori, for example as homogeneous sampling $t_k=kh$, for some step size $h>0$.

The time-triggered sampling operator $S_{\mathcal T}:C([0,\tau),\mathbb R^n)\to L^\infty_{\rm loc}([0,\tau),\mathbb R^n)$ is defined for a sampling sequence $\mathcal T=(t_k)_{k\in\mathbb N_0}$ with $t_0=0$, $t_k<t_{k+1}$ and $t_k\rightarrow\infty$ as $k\to\infty$ via
\[
  (S_{\mathcal T}u)(t)=u(t_k),\qquad \forall\,t\in[t_k,t_{k+1}).
\]
The respective closed-loop system is given by replacing $S_M$ with $S_{\mathcal T}$ in~\eqref{eq:System},~\eqref{eq:Controller}. The following theorem establishes a growth condition on the unknown nonlinear dynamics that is sufficient to conclude global existence of the solution for any sampling sequence $\mathcal{T}$.

\begin{thm}\label{thm:Fcondition}
Let $F$ be as in Step~3 of the proof of Theorem~\ref{Thm:event}, $\varphi(0)\Vert y^0-y_{\text{ref}}(0)+I^0\Vert <1$, and assume that the growth condition
\begin{equation}\label{eq:Fcondition}
\int_1^\infty\frac{dr}{1+r+F(r)}=\infty
\end{equation}
holds. Then, for every sampling sequence $\mathcal T$ with $t_k\to\infty$, there exists a global solution of the closed-loop system~\eqref{eq:System},~\eqref{eq:Controller} with $S_{\mathcal T}$ replacing $S_M$, and it satisfies $\varphi(t)\|e(t)\|<1$ for all $t\geq 0$. Moreover, $u$ and $I$ are bounded on every compact interval.
\end{thm}
\begin{proof}
The existence of a maximal solution $(y,I):[0,\omega)\to\mathbb R^n\times\mathbb R^n$ for some $\omega > 0$ can be shown by concatenation of the maximal solutions on the intervals $[t_k,\tilde{\omega})$ as shown in Step~1 of the previous proof. Let $R:[0,\omega) \rightarrow [0,\infty)$, $R(t) = \max_{s \leq t} \Vert I(s) \Vert$, and observe that by $\varphi(s)\Vert e(s) \Vert < 1$ on the solution of the closed-loop system, $$\Vert y(s) \Vert \leq \left\Vert \frac{1}{\varphi} \right\Vert_{\infty} + \Vert y_{\text{ref}} \Vert_{\infty} + R(t)$$ for all $s \leq t$, hence $\Vert f(y(s)) \Vert \leq F(R(t))$.

\emph{Step 1:} Derive a bound on $\varphi(s)\Vert e(s) \Vert$ and $u$ in terms of $R(s)$.
Define $$\lambda(r):=\left\|\frac{\dot\varphi}{\varphi}\right\|_\infty+\|\varphi\|_\infty(\|\dot y_{\text{ref}}\|_\infty+F(r)+\alpha r),$$ fix $\varepsilon_0 \in (0,1)$ with $\varphi(0)\|e(0)\|<\varepsilon_0$, and choose $\varepsilon(r) \in (0,1)$ such that $\frac{\varepsilon(r)^2}{1-\varepsilon(r)^2} = \max \left\{ \lambda(r), \frac{\varepsilon_0^2}{1-\varepsilon_0^2} \right\}$ for all $r \geq 0$. Then~\eqref{eq:relError} with $\tilde{M}$ replaced by $R(t)$ shows that $$\varphi(s)\|e(s)\|\le\varepsilon(R(t))$$ holds for all $s\le t$. 
Similarly, we can define $$ \Gamma(r) \coloneq \left\Vert \frac{1}{\varphi} \right\Vert_{\infty} \left( \lambda(r) + \frac{1}{1-\varepsilon_0^2}\right) $$ to obtain
\begin{equation*}
\Vert u(s)\Vert \leq \left\Vert \frac{1}{\varphi}\right\Vert_\infty\frac{\varepsilon(R(t))}{1-\varepsilon(R(t))^2}\leq \Gamma(R(t))
\end{equation*}
for all $s \leq t$. Observe that by the definition of $\lambda$, there exists $K >0$ such that $\Gamma(r) \leq K(1+r+F(r))$ holds for all $r \geq 0$.

\emph{Step 2:} Global existence of the solution.
Since $$\frac{d}{dt} \Vert I(t)\Vert \leq \Vert u(t)-S_{\mathcal{T}}u(t)\Vert \leq 2\Gamma(R(t))$$ we get
\begin{align*}
  R(t) = \max_{\tau \le t} \|I(\tau)\| &\leq \Vert I^0\Vert+2 \max_{\tau \le t} \int_0^\tau \Gamma(R(s)) ds\\
  &\le \Vert I^0\Vert+2\int_0^t \Gamma(R(s)) ds.
\end{align*}
Define $G:[\Vert I^0 \Vert, \infty) \rightarrow [0,\infty)$, $G(r) = \int_{\Vert I^0 \Vert}^r \frac{ds}{\Gamma(s)}$, and observe that~\eqref{eq:Fcondition} and the estimate $\Gamma(r) \leq K(1+r+F(r))$ show that $G$ is bijective. Since $\Gamma$ is nondecreasing, the Bihari-LaSalle inequality then yields that 
\[
  \forall\, t\ge 0:\ R(t) \leq G^{-1}(2t).
\]
Hence, $R$, and therefore $u$, are locally bounded. In particular, on each bounded interval $[0,T)$, there are only finitely many sampling times, so considering the last sampling interval, the solution remains within a compact subset of $\mathcal{D}$ defined in Step~1 of the proof of Theorem~\ref{Thm:event}, contradicting maximality for this sampling interval. Hence, the solution exists globally. 
\end{proof}

\begin{remark}
Condition~\eqref{eq:Fcondition} holds, for example if $F$ has a linear growth bound. The bounds on $u$ and $I$ are independent of the choice of the sampling sequence. Since there is no uniform bound on $\Vert I \Vert$, the tracking error $\Vert y(t) - y_{\text{ref}}(t) \Vert \leq \frac{1}{\varphi(t)} + G^{-1}(2t)$ may grow unboundedly as $t \rightarrow \infty$.
\end{remark}

\section{Simulations}

To illustrate the theoretical results, we consider the controlled Lorenz system
\begin{equation*}
\dot{y}(t) = f(y(t)) + S_Mu(t), \quad y(0) = y^0,
\end{equation*} 
where $f:\R^3 \rightarrow \R^3$ with 
\begin{equation*}
f(y_1,y_2,y_3) = \begin{pmatrix} \sigma(y_2 - y_1) \\ y_1(\rho - y_3) - y_2 \\ y_1y_2 - \beta y_3 \end{pmatrix}
\end{equation*}
with standard parameters $\sigma = 10, \rho = 28$ and $\beta = \frac{8}{3}$. Without control, this system exhibits chaotic behavior, making it a suitable stress test for the controller performance when aiming for reference tracking with performance constraints. The function $f$ is treated as unknown to the controller, and the funnel controller is designed without any knowledge of the system dynamics.

\begin{figure}[t]
\centering
\setlength{\abovecaptionskip}{0pt}
\setlength{\belowcaptionskip}{-10pt}
\includegraphics[width=\columnwidth]{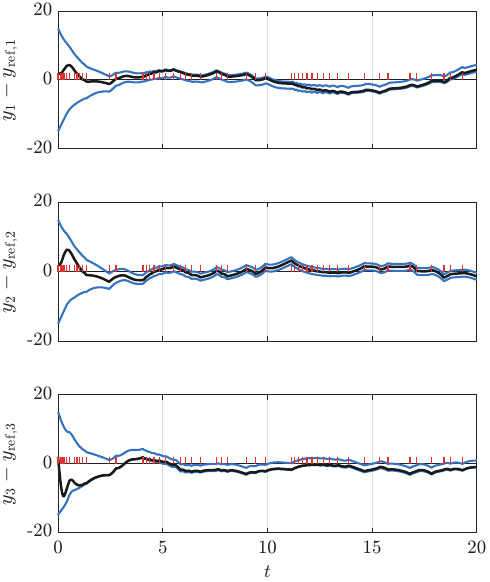}
\caption{Tracking error evolution within the shifted funnel boundaries of the closed-loop system with event-triggered funnel controller.}
\label{fig:Funnel}
\end{figure} 

\begin{figure}[t]
\centering
\setlength{\abovecaptionskip}{0pt}
\setlength{\belowcaptionskip}{-10pt}
\includegraphics[width=\columnwidth]{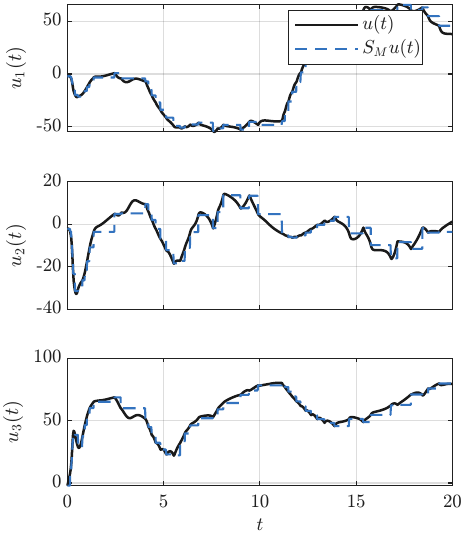}
\caption{Continuous computed signal $u(t)$ and applied zero-order hold input signal $S_Mu(t)$.}
\label{fig:Input}
\end{figure}

The control objective is to track the reference trajectory
\begin{equation*}
y_{\text{ref}}(t) = \begin{pmatrix} 3\sin(0.5t) \\ 7\sin(0.3t) \\ 4\sin(0.7t) + 25  \end{pmatrix},
\end{equation*}
with constraints described by the funnel function $\varphi(t) =  (14 e^{-t} + 1)^{-1} $. The initial conditions are $I^0 = 0 \in \R^3$ and $y^0 = (2,2,27)^{\top}$ satisfying the funnel condition initially. We choose the tuning parameters $\alpha = 2$ and $ M = 10$. All simulations have been conducted using MATLAB (solver: \texttt{ode45}, rel.tol.: $10^{-6}$, abs.tol.: $10^{-9}$, max.step: $0.01$).

\begin{figure}[t]
\centering
\setlength{\abovecaptionskip}{0pt}
\setlength{\belowcaptionskip}{-15pt}
\includegraphics[width=\columnwidth]{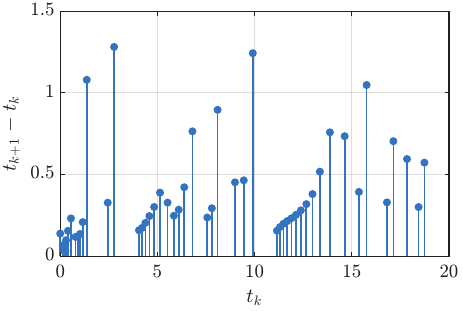}
\caption{Lengths $t_{k+1} - t_k$ of the sampling intervals generated by the event-triggered controller.}
\label{fig:Intervals}
\end{figure}

The evolution of the tracking error $y(t)-y_{\text{ref}}(t)$ within the funnel boundaries $[-I_j(t)-\varphi(t)^{-1}, -I_j(t) + \varphi(t)^{-1}]$ for $1\leq j\leq 3$ is depicted in Figure~\ref{fig:Funnel}. Note that the funnel constraint is imposed on the Euclidean norm of the error, hence the depicted boundaries are conservative for each component of the vector. The continuous time signal $u(t)$ and the sampled input $S_Mu(t)$ that is applied via zero-order hold are depicted in Figure~\ref{fig:Input} componentwise. 

\begin{figure}[t]
\centering
\setlength{\abovecaptionskip}{0pt}
\setlength{\belowcaptionskip}{-20pt}
\includegraphics[width=0.9\columnwidth]{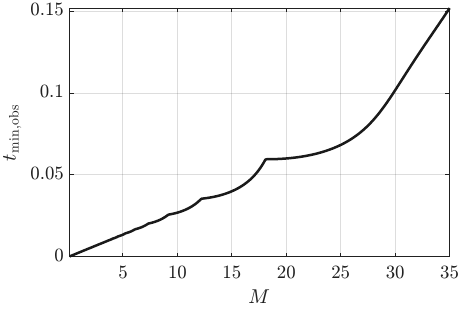}
\caption{Smallest observed interval length under variation of the tuning parameter $M$ while the other parameters remain fixed.}
\label{fig:Tuning}
\end{figure}

The lengths of the sampling intervals are shown in Figure~\ref{fig:Intervals}. In particular, the shortest observed sampling interval has a length of approximately $t_{\min, \text{obs}} \approx 0.027$, while the theoretical lower bound on the lengths of the sampling intervals can be determined to be approximately $t_{\min} \approx 6.8 \times 10^{-13}$. This shows the conservativeness of this theoretical bound. This gap reflects the worst-case nature of each estimate in the derivation of the theoretical bound $t_{\min}$.

Lastly, Figure~\ref{fig:Tuning} shows the smallest observed interval length $t_{\min,\text{obs}}$ in dependence of the tuning parameter $M$. While an increase in the tuning parameter $M$ causes an increase in the smallest observed interval length $t_{\min,\text{obs}}$ in this example, the theoretical bound derived in Theorem~\ref{Thm:event} is generally non-monotone in $\tilde{M}$.

\section{Conclusion}

In this paper, we have proposed an event-triggered zero-order hold implementation of funnel control for output reference tracking of unknown nonlinear systems by incorporating a dynamic auxiliary variable adjusting the reference trajectory. We have used the lower bound on the lengths of the sampling intervals to exclude Zeno behavior, and have derived a growth bound on the nonlinear dynamics that ensures the global existence of the solution for arbitrary, predefined sampling schedules. Future work may address the extension to systems of arbitrary relative degree or adaptive tuning strategies for the design parameters.

\small
\bibliographystyle{IEEEtran}
\bibliography{lit}

\end{document}